\documentclass[preprint,3p,12pt]{elsarticle}

\usepackage{color}
\usepackage{amssymb}
\usepackage{amsmath}
\usepackage{amsthm}

\usepackage{epsfig}
\usepackage{ulem}
\usepackage{epic}
\usepackage{tabls}
\usepackage{booktabs}
\usepackage{threeparttable}
\usepackage{pict2e}

\usepackage{float}
 \usepackage{subfig}

\usepackage{hyperref}

\hypersetup{backref,colorlinks=true,linkcolor=blue}

\usepackage[vlined]{algorithm2e}

\begin{document}
\begin{frontmatter}



\title{Reduced-Order Models for Thermal   Radiative Transfer \\
  Based on  POD-Galerkin Method and Low-Order Quasidiffusion Equations }

\author[ncsu,ncsu1]{Joseph M. Coale}
\author[ncsu,ncsu2]{Dmitriy Y. Anistratov}
\address[ncsu]{Department of Nuclear Engineering
North Carolina State University Raleigh, NC}
\address[ncsu1]{jmcoale@ncsu.edu}
\address[ncsu2]{anistratov@ncsu.edu}

\begin{abstract}
This paper presents  a new technique for developing reduced-order models (ROMs) for
nonlinear  radiative transfer problems in high-energy density physics.
 The  proper orthogonal decomposition  (POD) of photon intensities  is
  applied
  to obtain global basis functions for the Galerkin projection  (POD-Galerkin)
 of the time-dependent multigroup Boltzmann transport equation (BTE) for photons.
  The POD-Galerkin solution of the BTE is used to determine the quasidiffusion (Eddington) factors
  that yield  closures  for the nonlinear system of  (i) multilevel low-order quasidiffusion (VEF) equations
  and (ii) material energy balance equation.
Numerical results are presented  to demonstrate accuracy of the ROMs  obtained with different low-rank approximations of intensities.
 \end{abstract}


\end{frontmatter}

\section{Introduction}

In this paper, we develop  reduced order models (ROMs) for  the basic thermal radiative transfer (TRT) problem that neglects material motion, scattering and heat conduction.
Problems in 1D slab geometry are considered.
The TRT problem is defined by the time-dependent multigroup Boltzmann transport  equation (BTE) given by
\begin{gather}
	\frac{1}{c} \frac{\partial I_g  }{\partial t}(x, \mu, t)
	+ \mu \frac{\partial  I_g}{\partial  x} (x, \mu, t)
	+ \varkappa_{g}(T)I_g(x, \mu, t) =
	2\pi \varkappa_{g}(T) B_g(T) \label{bte}\\
	x \in [0,X], \quad \mu \in [-1,1], \quad g=1,\ldots,N_g, \quad  t \in [0,  t_{end}], \nonumber\\
	I_g|_{{\mu>0}\atop{x=0}} =  I_g^{\text{in}+}, \quad I_g|_{{\mu<0}\atop{x=X}} =  I_g^{\text{in}-}, \quad I_g|_{t=0} =  I_g^0,
\end{gather}
and the  material energy balance (MEB) equation
\begin{equation}\label{meb}
\frac{\partial \varepsilon(T)}{\partial t} =
\sum_{g=1}^{N_g}  \varkappa_{g}(T) \Big(\int_{-1}^1 I_g(x, \mu, t)d  \mu  -  4\pi B_g(T) \Big)   \, ,
\quad \left. T \right|_{t=0} =  T_0.
\end{equation}
Here
$I_g$ is the group specific photon intensity;
$x$ is the spatial position;
$\mu$ is   the direction cosine of particle motion;
$g$ is the index of photon frequency group;
$N_g$ is the number of frequency groups;
$t$ is time;
$\varkappa_g$ is the group material opacity;
$T$ is the material temperature;
$\varepsilon$  is the material energy density and
 $B_g$ is the group Planck black-body distribution  function.
The BTE  describes propagation and absorption of photons in matter, and emission of photons with the black-body spectrum.
The MEB equation models change in material energy as the result of absorption and emission of photons.

Particle transport problems have high dimensionality; Discretization of the BTE in the phase space and time results in a problem with a very large number of degrees of freedom (DoF).
This has stimulated active research on the development of ROMs for the BTE and its related class of problems
\cite{jc-dya-m&c2019,jc-dya-tans2019,ragusa-morel-m&c2019,ryan-jcp2020,Choi-2021}.
In this study, the new ROMs for TRT problems  \eqref{bte}-\eqref{meb} are based
on the proper orthogonal decomposition (POD) and projection approach  \cite{sirovich-1987,berkooz-holmes-lumley-1993,volkwein-2002}.
The POD creates an optimal basis to represent  dynamics of a   system
based on a set of  collected data \cite{Benner-siam-2015}.
Specifically, the reduced basis for a Galerkin expansion of intensities over the phase space is  generated by the POD of a collection of
vectors of numerical transport solutions over time intervals of the TRT problem.
  The BTE is projected onto this basis.
The   POD-Galerkin projected BTE can be viewed as a discretization scheme based
on a set of problem-specific global basis functions.
The  projected BTE is then coupled with the multilevel
nonlinear system of governing moment equations consisting of (i) low-order quasidiffusion (aka VEF) equations for the group and total radiation energy densities and fluxes and (ii) MEB equation.
The moment equations are derived by a nonlinear projection of the BTE (Eq. \eqref{bte})
using exact closures by means of  quasidiffusion (Eddington) factors.
The POD-Galerkin expansion of intensities is used to compute the quasidiffusion (QD) factors.

The reminder of the paper is organized as follows.
 In Sec. \ref{sec:POD-G} POD-Galerkin projection of the BTE is formulated.
 In Sec. \ref{sec:mlqd} the ROM based on multilevel low-order QD equations is described.
 Numerical results are presented in Sec. \ref{sec:num}.
 We conclude with a brief discussion in Sec. \ref{sec:end}.

\section{\label{sec:rom} Reduced-Order Model for TRT}

\subsection{\label{sec:POD-G} POD-Galerkin Projection of BTE}

To discretize the BTE  \eqref{bte} we apply
 (i) the method of discrete ordinates (MDO) for the angular variable,
 (ii) the backward Euler (BE)  scheme  for time integration,
 and  (iii)   the simple corner-balance (SCB) method
   for  approximation in space \cite{mla-ttsp-1997} to obtain
\begin{equation}\label{bte-disc}
 \frac{1}{c\Delta t^n}  \big( \mathbf{I}^n - \mathbf{I}^{n-1} \big)
+ \mathcal{L}_h \mathbf{I}^n
+   \mathcal{K}_{h}^n(T) \mathbf{I}^n
=   \mathbf{Q}^n(T) \, ,
\end{equation}
where $n$ is the time step index,
$\mathbf{I}^n=((\boldsymbol{I}^n_1)^\top\ \dots\ (\boldsymbol{I}^n_{N_g})^\top)^\top\in \mathbb{R}^D$
is the solution vector at $t=t^n$,
$\boldsymbol{I}^n_g \in \mathbb{R}^{2N_x N_{\mu}}$ is the vector of group intensities,
$ D=2N_x N_{\mu}N_g $ is the number of DoF in the phase space,
$N_x$  is the number of  spatial mesh cells,
$N_{\mu}$  is the number of discrete  angular directions,
$\mathcal{L}_h$ and $\mathcal{K}_{h}$  are the discrete operators that define
 approximation by the MDO and the SCB  scheme,
 $\mathbf{Q}^n(T)$  is the  vector of the right-hand side,
 $\Delta t^n$ is the $n^\text{th}$ time step.

The numerical solution of the discretized BTE \eqref{bte-disc}
on a given phase-space grid is used to form
  a database matrix $\mathbf{A} =  [\mathbf{I}^1, \ldots,\mathbf{I}^{N_t}] \in \mathbb{R}^{D\times N_t}$ of the solution snapshots  computed over $N_t$ time steps.
The database $\mathbf{A}$ is used to form a POD basis $\{\boldsymbol{u}_\ell\}_{\ell=1}^r$
 with  $r \ll D$ that gives an optimal approximation of $\mathbf{A}$  and solves the following optimization problem \cite{volkwein-2013}:
\begin{equation}
	\min_{\boldsymbol{u}_1,\dots,\boldsymbol{u}_r}\sum_{n=1}^{N_t} \Delta t^n \bigg\| \mathbf{I}^n - \sum_{\ell=1}^r\langle\mathbf{I}^n,\boldsymbol{u}_\ell\rangle_W\boldsymbol{u}_\ell \bigg\|_W^2 \, , \label{pod_opt}
\end{equation}
where the norm $\|\cdot\|_W^2$ is defined by the spatial and angular discretization of the BTE.
The weighted inner product specific  to the SCB and MDO discretization is given by  $\big<  \boldsymbol{u}_{\ell},\boldsymbol{u}_{\ell'} \big>_W = \boldsymbol{u}_{\ell}^\top \mathbf{W} \boldsymbol{u}_{\ell'}$ with
\begin{equation}
	\mathbf{W} = \bigoplus_{g=1}^{N_g}\bigoplus_{m=1}^{N_\mu} w_m\hat{\mathbf{W}}^x, \qquad \mathbf{W}\in\mathbb{R}^{D \times D} \, .
\end{equation}
Here $w_m$ are the angular quadrature weights,  and $\hat{\mathbf{W}}^x = \bigoplus_{i=1}^{N_x} \frac{\Delta x_i}{2} \mathbb{I}$, where $\mathbb{I} \in \mathbb{R}^{2\times2}$ is the identity matrix. We form the   weighted data matrix
\begin{equation}
\mathbf{\hat A} =\mathbf{W}^{1/2}\mathbf{A}\mathbf{D}^{1/2}\,  \quad \mbox{with} \quad
\mathbf{D}=\text{diag}(\Delta t^1,\dots,\Delta t^{N_t})
\end{equation}
 and find its singular value decomposition (SVD) to get
\begin{equation}
\mathbf{\hat A} = \hat{\mathbf{U}}   \hat{\mathbf{\Sigma}}   \hat{\mathbf{V}}^\top \, ,
\end{equation}
 where
 $\hat{\mathbf{U}}= [\hat{\boldsymbol{u}}_{1}, \dots,\hat{\boldsymbol{u}}_{d}]\in\mathbb{R}^{D \times d}$  holds the left singular vectors of $\mathbf{\hat A}$ in its columns,
 $\hat{\mathbf{V}}  \in\mathbb{R}^{N_t \times d}$ is the matrix of the right singular vectors,
$\hat{\mathbf{\Sigma}}=\text{diag}( \sigma_1,\dots, \sigma_d) \in \mathbb{R}^{d \times d}$
is the matrix of singular values, $d=\min(D,N_t)$ is the  rank of $\mathbf{\hat A}$.
The POD basis  $\mathbf{U} = [\boldsymbol{u}_{1}, \dots,\boldsymbol{u}_{d}] \in \mathbb{R}^{D \times d}$ satisfying Eq. \eqref{pod_opt} is given by
 \begin{equation} \label{podbasis}
 \mathbf{U}  = \mathbf{W}^{-1/2} \hat{\mathbf{U}} \, .
 \end{equation}
We now formulate Galerkin ansatz  expanding the intensities in the POD basis  \eqref{podbasis} \cite{volkwein-2002,volkwein-2013}
\begin{equation} \label{I-podg}
 \mathbf{I}^{u}_r(t^n) = \sum_{\ell=1}^{r}  \lambda_{\ell}^n \boldsymbol{u}_{\ell}, \quad r\leq d .
\end{equation}
The discretized transport equation \eqref{bte-disc} is projected onto the POD basis
  to derive  the POD-Galerkin (POD-G) projected BTE ($\ell=1,\dots,r$) given by
 \begin{equation} \label{podg-bte}
 \frac{1}{c\Delta t^n}       \big( \lambda_{\ell}^n   -  \lambda_{\ell}^{n-1} \big)
  +   \sum_{\ell'=1}^{r}    \lambda_{\ell'}^n  \big<  \boldsymbol{u}_{\ell}, \mathcal{L}_h   \boldsymbol{u}_{\ell'}\big>_W
  +    \sum_{\ell'=1}^{r}   \lambda_{\ell'}^n  \big<  \boldsymbol{u}_{\ell}, \mathcal{K}_{h}^n(T)   \boldsymbol{u}_{\ell'} \big>_W =
  \big<  \boldsymbol{u}_{\ell},\mathbf{Q}^n(T) \big>_W \, ,
\end{equation}
where  it is taken into account that $\big<  \boldsymbol{u}_{\ell'},\boldsymbol{u}_{\ell} \big>_W=\delta_{\ell,\ell'}$.
The POD basis is a global one. This yields a non-sparse  system of equations for  the coefficients $\lambda_{\ell}^n$ at $t=t^n$. However, it will be shown below that the ROMs based on the POD-G  projected BTE \eqref{podg-bte} are accurate for $r  \ll D$.

\subsection{\label{sec:mlqd}  ROM Based on Multilevel Low-Order QD Equations}

The  system of  equations of the multilevel QD (MLQD) method for the TRT problem  \eqref{bte}-\eqref{meb} is  derived by  a nonlinear projection  in angular and frequency (photon energy) spaces.
It    is defined  by  the following  low-order equations \cite{gol'din-1964,gol'din-1972,PASE-1986}:
\begin{enumerate}
\item   The multigroup  low-order QD (LOQD)  equations  for the angular moments  given by
\begin{subequations}\label{mloqd}
\begin{gather}
 \frac{\partial E_g }{\partial t} +  \frac{\partial F_g }{\partial x} +
 c\varkappa_{g}(T) E_g    =   4\pi \varkappa_{g}(T) B_g(T) \, ,\\
 \frac{1}{c} \frac{\partial  F_g}{\partial t} +
  c  \frac{\partial  }{\partial x} \Big( f_g[I]    E_g \Big) + \varkappa_{g} (T)  F_g
   =   0 \, ,
\end{gather}
\end{subequations}
where
$E_g = \frac{1}{c} \int_{-1}^1 I_g d \mu$
is the   group radiation energy density,
$F_g =  \int_{-1}^1 \mu I_g d \mu$ is the group radiation flux, and
\vspace{-0.15 cm}
\begin{equation} \label{f}
	f_g\big[ I \big]  = \left.\int_{-1}^1 \mu^2   I_g  d \mu \right/  \int_{-1}^1   I_g  d \mu
\end{equation}
is the group QD (Eddington) factor that provides  closure of the BTE and multigroup LOQD equations.
This closure is exact when $f_g$ is defined by the solution of the BTE \eqref{bte} according to the QD (VEF) method \cite{gol'din-1964,auer-mihalas-1970}.
\item   The effective grey LOQD equations for
the  total radiation energy density $E=  \sum_{g=1}^{N_g} E_g$ and
the  total   flux $F=   \sum_{g=1}^{N_g} F_g$  are given by
\begin{subequations}\label{gloqd}
\begin{gather}
	\frac{\partial E   }{\partial t}
	+ \frac{\partial F   }{\partial x}
	+ c \bar \varkappa_E  E
	= c  \bar \varkappa_B a_R T^4 \,  ,\label{gloqd-1}\\
	\frac{1}{c}\frac{\partial F   }{\partial t}
	+ c\frac{\partial (\bar f \big[I \big] E)   }{\partial x}
	+ \bar \varkappa_R  F + \bar \eta E = 0 \, , \label{gloqd-2}
\end{gather}
\end{subequations}
where the grey coefficients are
\begin{gather*}
	\bar{\varkappa}_{E}= \frac{\sum_{g=1}^{N_g} \varkappa_{g}E_g}
	{\sum_{g=1}^{N_g} E_g}   \, ,
	\quad
	\bar{\varkappa}_{B}= \frac{\sum_{g=1}^{N_g} \varkappa_{g} B_g}
	{\sum_{g=1}^{N_g} B_g}  \, ,
	\quad
	\bar{\varkappa}_{R}= \frac{\sum_{g=1}^{N_g} \varkappa_{g}|F_g|}
	{\sum_{g=1}^{N_g} |F_g|} \ ,\\
	\bar f = \frac{\sum_{g=1}^{N_g} f_g E_g}{\sum_{g=1}^{N_g} E_g}  \, ,
	\quad
   	\bar \eta  = \frac{\sum_{g=1}^{N_g} (\varkappa_{g} - \bar{\varkappa}_{R}) F_g}{\sum_{g=1}^{N_g} E_g} \, .
\end{gather*}
\end{enumerate}
The grey LOQD equations are coupled  with the MEB equation that is cast in grey form for the total energy density
\begin{equation}\label{eb-grey}
\frac{\partial \varepsilon(T)}{\partial t} =  c  \big(\bar \varkappa_E E - \bar\varkappa_B a_R T^4\big) \, .
\end{equation}

The   new ROM for TRT   combines the POD-G projected  BTE
with nonlinear projection in angular variable and photon energy
via the hierarchy  of  low-order QD equations for moments of the intensity.
It is defined by the following set of equations:
\begin{itemize}
\item the POD-G projected BTE (Eq. \eqref{podg-bte}) the solution of which gives compressed  representation
of the intensities in the phase space,
\item the multigroup   LOQD equations (Eq. \eqref{mloqd}),
 where the  QD factors are defined by  the POD-G   expansion of intensities of rank $r$ and hence
\begin{equation}
 f_g^{u} = f_g\big[ \mathbf{I}^{u}_r \big] \, ,
\end{equation}
\item the effective grey LOQD equations (Eq. \eqref{gloqd}) and the MEB equation in the grey form (Eq. \eqref{eb-grey}).
\end{itemize}
The QD factor $f_g^{u}$  defines  an approximate closure for the group LOQD  equations
providing further data compression of intensities and the next level of  reduction of dimensionality for the TRT problem.
Hereafter  we refer to this ROM as the QD-PODG model, whose iterative algorithm for solving TRT problems is outlined in Algorithm \ref{alg:rom_alg}.
Temporal discretization of the  LOQD and MEB equations  (Eqs.   \eqref{mloqd}, \eqref{gloqd}, and \eqref{eb-grey}) is based on  the  BE  time integration method.
 The multigroup LOQD
 equations  are discretized in space by means of a second-order finite volume (FV) method \cite{dya-jcp-2019}.
 The spatial discretization of the grey LOQD equations   is algebraically  consistent with  the  discretized multigroup LOQD equations.

The coefficients of the POD-G projected BTE   explicitly depend on $T$  through group opacities and the Planckian emission term.
This makes Eq. \eqref{podg-bte} an integral part of the nonlinear multilevel system of LOQD equations by means of which they are coupled to the MEB equation.
This feature allows these equations to be used in the development of parameterized ROMs for TRT.
One can generate the POD-G basis $\{\boldsymbol{u}_\ell\}_{\ell=1}^r$ for a base case TRT problem and
use the QD-PODG model with this basis to solve TRT problems with different parameters, for example,
a perturbed  spectrum of incoming radiation.

\begin{algorithm}[ht!]
\DontPrintSemicolon
	\SetAlgoLined
	\While{$t_n \leq t^\text{end}$}{\vspace*{.1cm}
		$n = n+1$\\
		$T^{(0)} = T^{n-1}$\\
		\While{ $\|T^{(k)} - T^{(k-1)}\| > \epsilon_T\|T^{(k)}\|, \ \ \|E^{(k)} - E^{(k-1)}\| > \epsilon_E\|E^{(k)}\|$ }{
			$k=k+1$\;
            Update $\varkappa_g,B_g$  using $T^{(k-1)}$\;
			Solve   POD-G projected BTE (Eq. \eqref{podg-bte})   given $\{\boldsymbol{u}_\ell\}_{\ell=1}^r, \ T^{(k-1)}$ to compute $\{\lambda_\ell^{(k)}\}_{\ell=1}^r$\;
Compute $\mathbf{I}^{u (k)}_r =\sum_{\ell=1}^{r}  \lambda_{\ell}^{(k)} \boldsymbol{u}_{\ell}$\;
			Compute  $f_g^{u(k)}= f_g\big[\mathbf{I}^{u(k)}_r \big]$\;
			\While{ $\|T^{(k,s)} - T^{(k,s-1)}\| > \epsilon_T\|T^{(k,s)}\|, \ \ \|E^{(k,s)} - E^{(k,s-1)}\| > \epsilon_E\|E^{(k,s)}\|$ }{
				$s=s+1$\;
                Update $\varkappa_g,B_g$  using $T^{(k,s-1)}$ \;
				Solve Eqs. \eqref{mloqd}   given $f_g^{u,(k)}$ to compute $E_g^{(k,s)}, F_g^{(k,s)}$ \;
				Compute grey coefficients $\bar{\varkappa}_{E}^{(k,s)}, \bar{\varkappa}_{B}^{(k,s)}, \bar{\varkappa}_{R}^{(k,s)}, \bar f^{(k,s)}, \bar \eta^{(k,s)}$\;
				Solve Eqs. \eqref{gloqd} and \eqref{eb-grey} to compute  $E^{(k,s)}, F^{(k,s)}, T^{(k,s)}$ \;
			}
		$T^{(k)}\leftarrow T^{(k,s)}$\;
		}
		$T^n\leftarrow T^{(k)}$, $\lambda_{\ell}^{n}  \leftarrow \lambda_{\ell}^{(k)}$\;
	}
	\caption{The algorithm of the QD-PODG model for solving TRT problems \label{alg:rom_alg}}
\end{algorithm}

\section{\label{sec:num} Numerical Results}

To analyze  the accuracy of the
 QD-PODG model, we use
the  problem based on the well-known Fleck-Cummings (F-C) test \cite{fleck-1971}.
A 1D slab of one material is defined as 6 cm thick ($X=6$).
The material  spectral opacity is   given by
$\varkappa_{\nu} = \frac{27}{(h\nu)^3}\big(1-e^{-\frac{h\nu}{kT}}\big)$.
The left boundary has incoming radiation with black-body spectrum $B_{\nu}$ at temperature
$kT_{in}=1$ keV
 and the right boundary is vacuum. The initial temperature of the slab is  $kT_0=1$ eV
  and the initial radiation distribution is given by the black-body spectrum at $T_0$.
The material energy density is a linear function of temperature $\varepsilon = c_\nu T$, where  $c_\nu= 0.5917a_RT_{in}^3$.
 The time interval of the problem  is $0 \le t  \le 6$ ns.
 A uniform time step is used $\Delta t~= 2 \times 10^{-2}$~ns and hence there are 300 time steps ($N_t=300$).
   The spatial mesh consists of a uniform $N_x=60$ cells with width $\Delta x = 0.1$ cm.
The angular mesh has  8  discrete directions ($N_\mu=8$).
The double $S_4$ Gauss-Legendre  quadrature set is used.
We define $N_g=17$ energy groups.
The parameters of convergence criteria for temperature and energy density are  $\epsilon_T=\epsilon_E=10^{-12}$, respectively.

The full-order model (FOM) for this TRT problem is   formulated as   the MLQD method
where the BTE and low-order QD equations are discretized as described above on the given grid in phase space and time.
 The number of  DoF   of  $\mathbf{I}^n$ at each instant of time $t^n$  is  $D=1.632\times10^{4}$.
The number of DoF in the phase space and time  for this FOM  is equal to $DN_t= 4.896\times10^{6}$.
The solution  to  the F-C test evolves in three distinct temporal stages: (i) rapid wave formation, (ii) wave propagation, and (iii) slow continual heating of the domain to steady state.
 A separate database is constructed by the FOM for each of these stages, whose temporal ranges are the following:
$t \in [0,0.3 \, \text{ns}]$  for $i=1$,
$t \in (0.3, 1.2 \,  \text{ns}]$  for $i=2$,
$t \in (1.2, 6 \, \text{ns}]$  for $i=3$.

\begin{figure}[t]
	\centering \hspace*{-.5cm}
	\subfloat[$\mathbf{A}_{1}\ (0\leq t\leq 0.3 \, \text{ns})$]{\includegraphics[width=.35\textwidth]{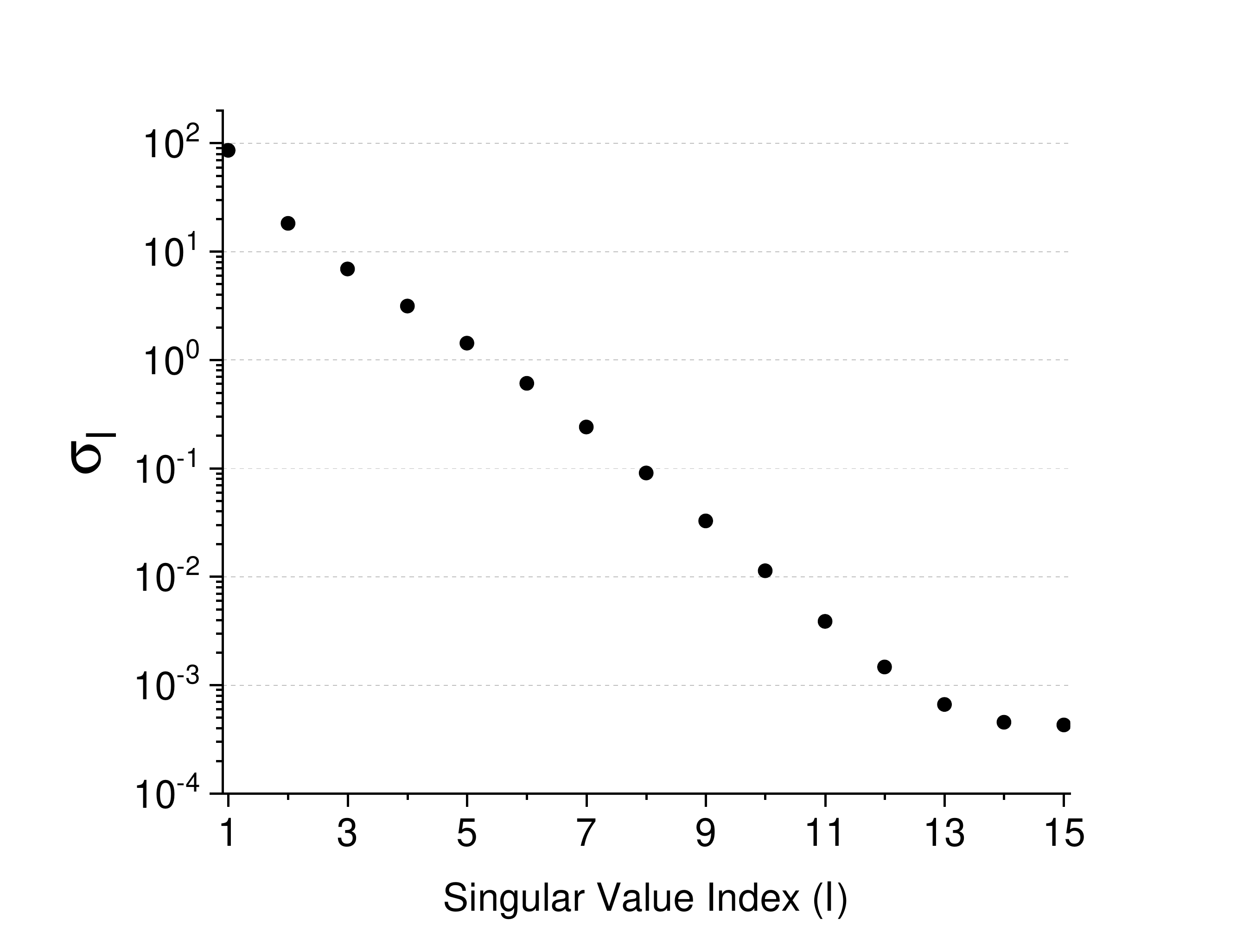}}
	\subfloat[$\mathbf{A}_{2}\ (0.3< t\leq 1.2 \, \text{ns})$]{\includegraphics[width=.35\textwidth]{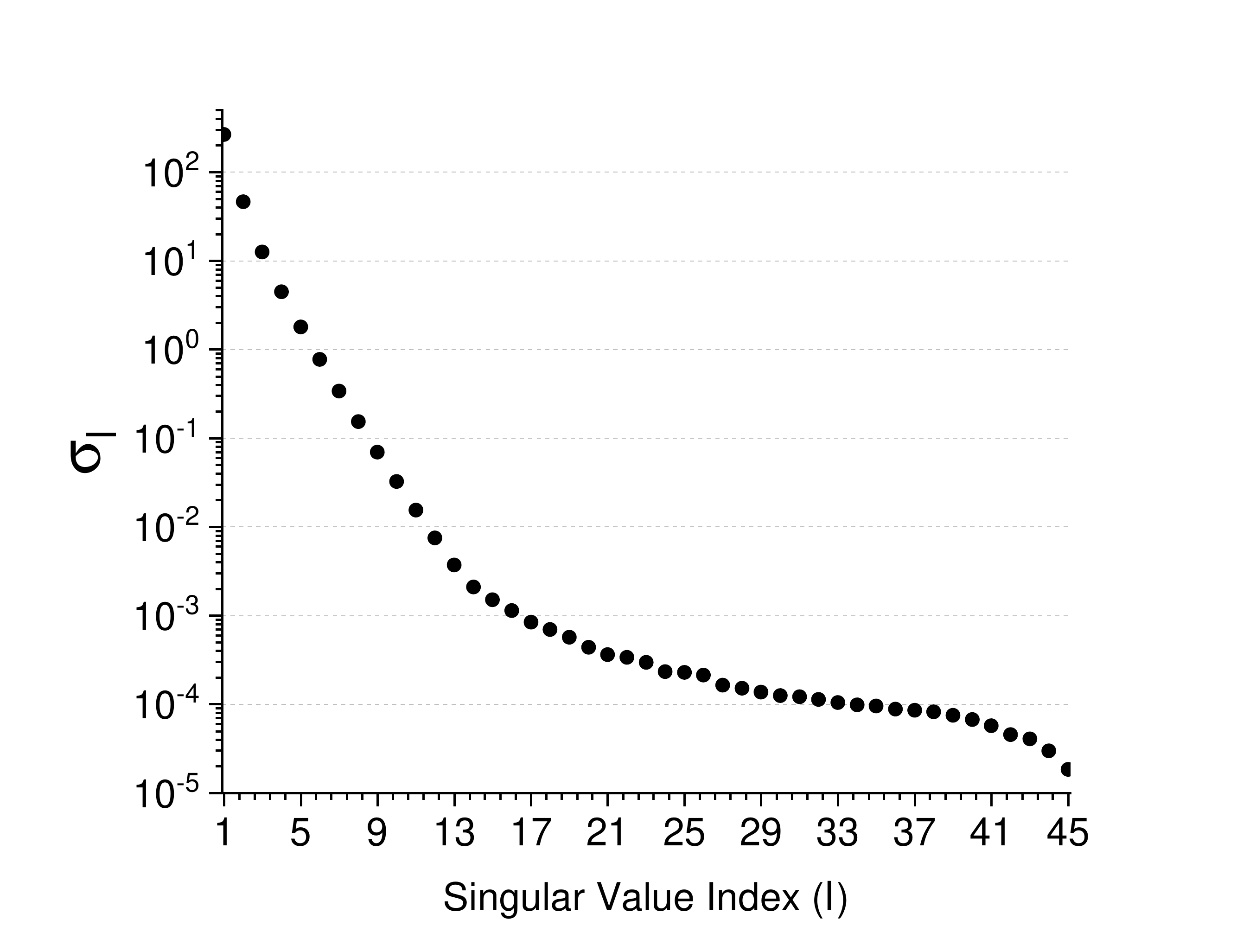}}
	\subfloat[$\mathbf{A}_{3}\ (1.2< t\leq 6 \, \text{ns})$]{\includegraphics[width=.35\textwidth]{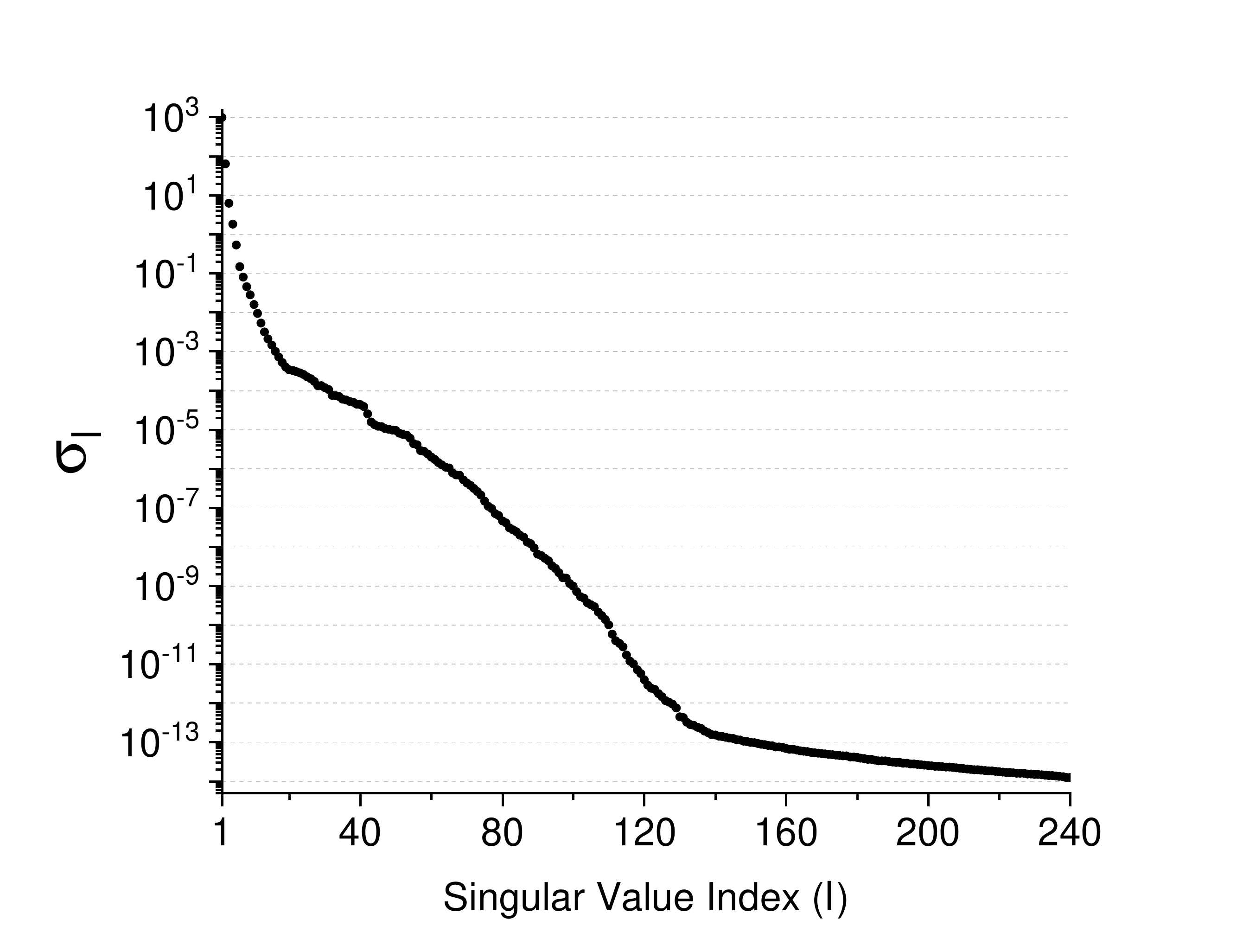}}
	\caption{Singular values of the database matrices of intensities
over three time subintervals of the problem \label{fig:svals}}
\end{figure}

The resulting database matrices that hold the set of discrete intensities for each of the three stages of the F-C test we denote by $\mathbf{A}_i   \in \mathbb{R}^{D \times N_{t,i}}$, $i=1,2,3$.
  The columns of each database are snapshots of the solution at $N_{t,1},N_{t,2},N_{t,3}$ instants of time, respectively,  ordered chronologically.
The full ranks $d_i$ of $\mathbf{A}_i$  are equal to $d_1=N_{t,1}=15,\ d_2=N_{t,2}=45,\ d_3=N_{t,3}=240$ respectively. The singular values $(\sigma_{\ell})$ of each of the three databases are depicted in Figure \ref{fig:svals}.
The first database shows a slow rate of decrease in magnitude of its singular values over the entire range,
 whereas the singular values of
the other two databases first experience rapid decrease  followed by  a plateau
 where the change in their magnitudes   slows significantly.
From the matrices $\mathbf{A}_i $, POD bases  $\{\boldsymbol{u}_{i,\ell}\}_{\ell=1}^{r_i}$, $i=1,2,3$ are calculated for each of these time intervals.

\begin{figure}[t]
	\centering
	\includegraphics[scale=0.3]{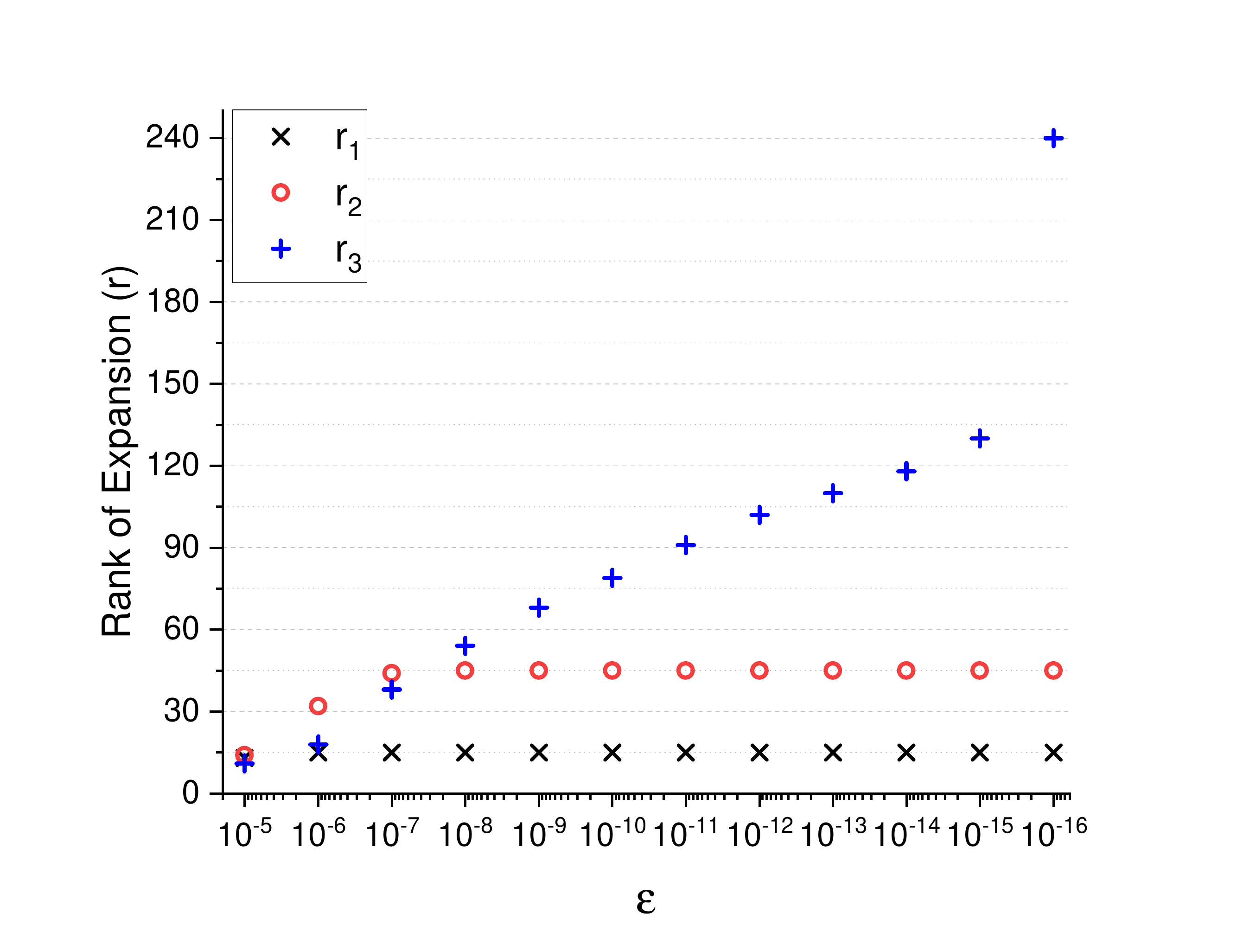}
	\caption{Rank of expansion for each database corresponding to different $\varepsilon$ \label{fig:exp_ranks}}
	\vspace*{-0.5cm}
	\centering
	\subfloat[Temperature]{\includegraphics[width=.48\textwidth]{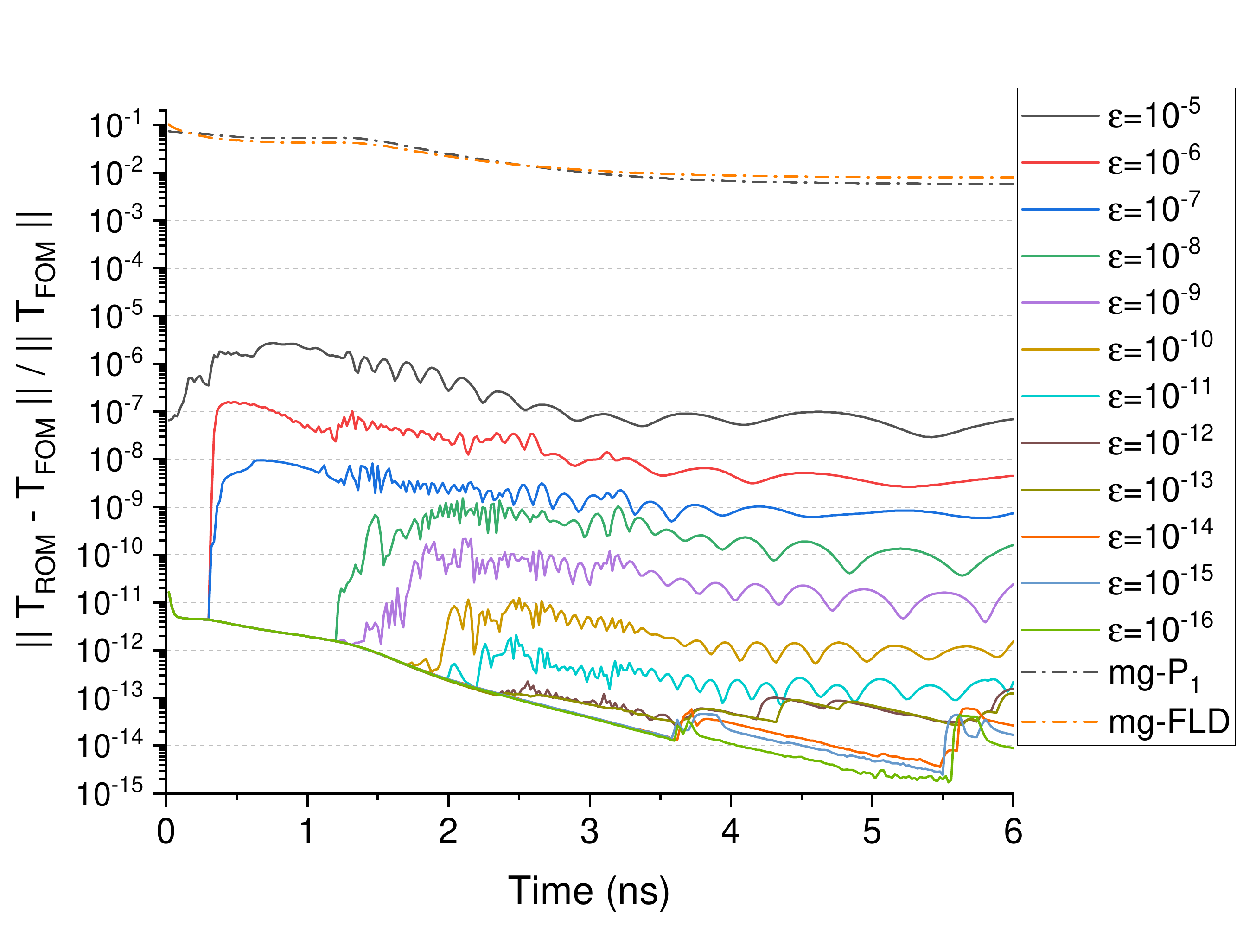}}
\hspace{0.2cm}
	\subfloat[Radiation Energy Density]{\includegraphics[width=.48\textwidth]{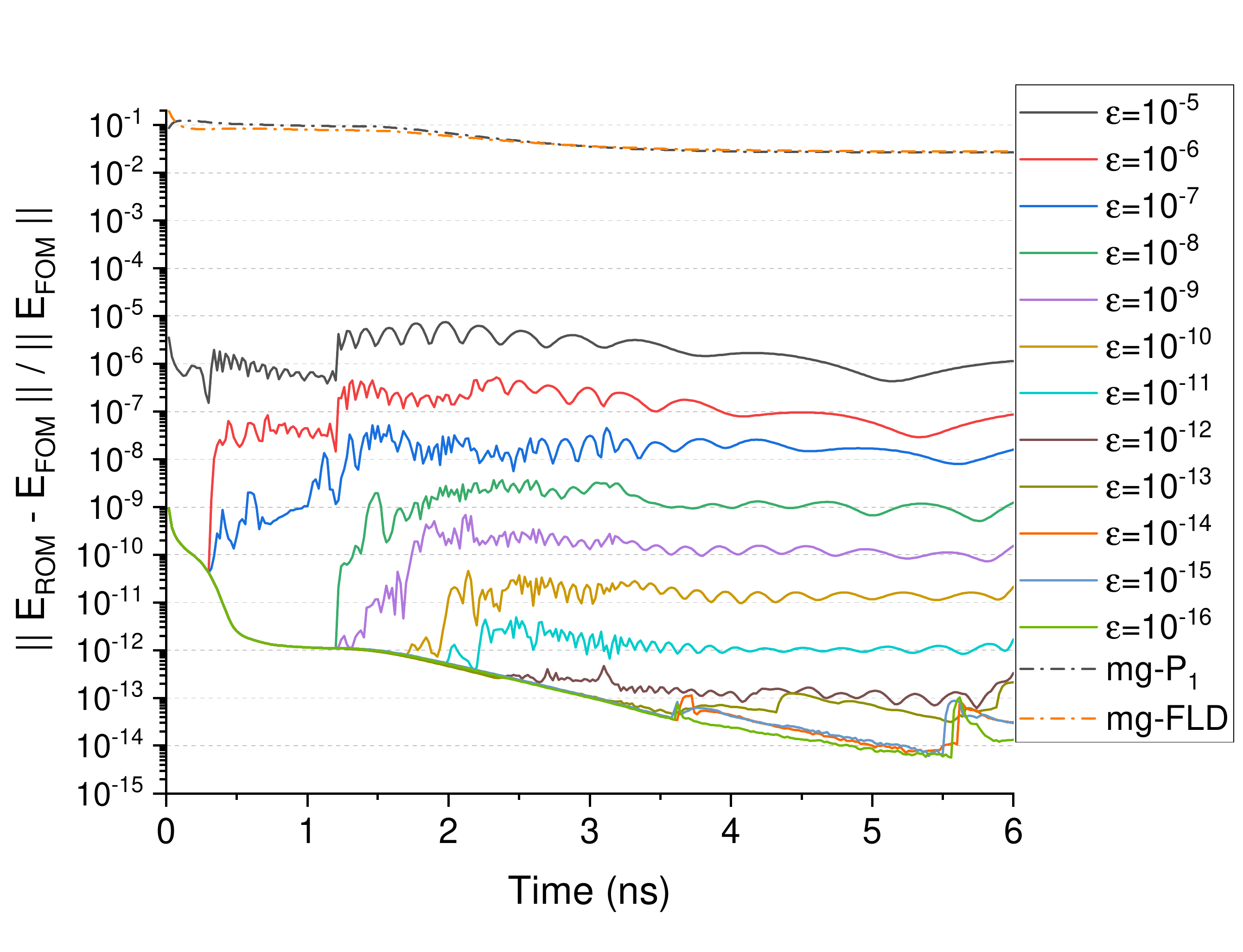}}
	\caption{Relative errors of the QD-PODG model mg-$P_1$ and mgFLD ROMs compared to the FOM solution in the 2-norm vs. time \label{fig:time_err}}
\end{figure}

We now solve the F-C test   with the QD-PODG model by expanding $I$ with each of our three POD bases based on the time frames they were generated for (e.g. we expand with $\{\boldsymbol{u}_{1,\ell}\}_{\ell=1}^{r_1}$ while $0\leq t \leq 0.3$ ns).
The ranks  $(r_1,r_2,r_3)$ of the expansion \eqref{I-podg}  are determined as the values  that satisfy the following criterion \cite{sirovich-1987}
\begin{equation} \label{eta-g}
	\bigg( \sum_{\ell=r_i+1}^{d_i} \sigma_{\ell}^2 \Big/
	\sum_{\ell=1}^{d_i}\sigma_{\ell}^2 \bigg)^{\frac{1}{2}}    < \varepsilon, \quad
\mbox{for} \ A_i \, , \ i=1,2,3,
\end{equation}
given some desired $\varepsilon$. The ranks found for $\varepsilon\in[10^{-5},10^{-16}]$ are shown in Figure \ref{fig:exp_ranks}.
The POD bases for $\mathbf{A}_1$ and $\mathbf{A}_2$ reach full rank $(r_1=15,\ r_2=45)$ at $\varepsilon= 10^{-6}$ and $\varepsilon= 10^{-8}$, respectively.
 Full-rank is not found for the basis of $\mathbf{A}_3\ (r_3=240)$ until $\varepsilon= 10^{-16}$.
This behavior is expected since compared to $\mathbf{A}_3$, the full ranks of $\mathbf{A}_1$ and $\mathbf{A}_2$ are relatively small. The singular values of both $\mathbf{A}_1$ and $\mathbf{A}_2$ also occupy a smaller range than for $\mathbf{A}_3$.
Another notable behavior is that $r_3<r_2$ for $\varepsilon< 10^{-8}$, indicating that the solution contained in the time range over which $\mathbf{A}_2$ was generated is the most difficult to represent with few POD modes.
This is to be expected given that $\mathbf{A}_2$ accounts for the solution during propagation of the radiation wave from the left boundary to the right, which is known to be a difficult phenomena for the POD to represent with low rank \cite{rowley-mardsen-2000,reiss-2004}. Let us note here that the rank of expansion for each timeframe in the F-C test is exactly the size of the linear system that   solves for the coefficients $\lambda_\ell$ (Eq. \eqref{podg-bte}). This means that when using $\varepsilon = 10^{-5}$ for instance, the largest linear system to solve in place of the BTE is a dense $r \times r$ system with $r=14$, which is of significantly lower dimensionality than the original BTE ($r=14 \ll D=1.632\times10^{4}$).

The errors of the QD-PODG model relative to the  FOM  solution on the F-C test in the 2-norm are displayed in Figure \ref{fig:time_err} for material temperature and radiation energy density vs. time. Each unique curve shows the relative error of
the ROM solution for a specific value of $\varepsilon$, ranging from $10^{-5}$ to $10^{-16}$.
 Note that we use the FOM solution as the reference  to compute errors against,  to determine how the ROM solution converges  to its training data. The MLQD discrete solution will converge to the multigroup TRT solution in the limit $N_x,N_\mu,N_t\rightarrow\infty$ and so we can postulate that if the solution of the QD-PODG model converges to the discrete FOM solution then it will too converge to the continuous solution given a database generated on a fine-enough grid.

\begin{figure}[t]
	\vspace*{-1cm}
	\centering
	\subfloat[Temperature]{\includegraphics[width=.5\textwidth]{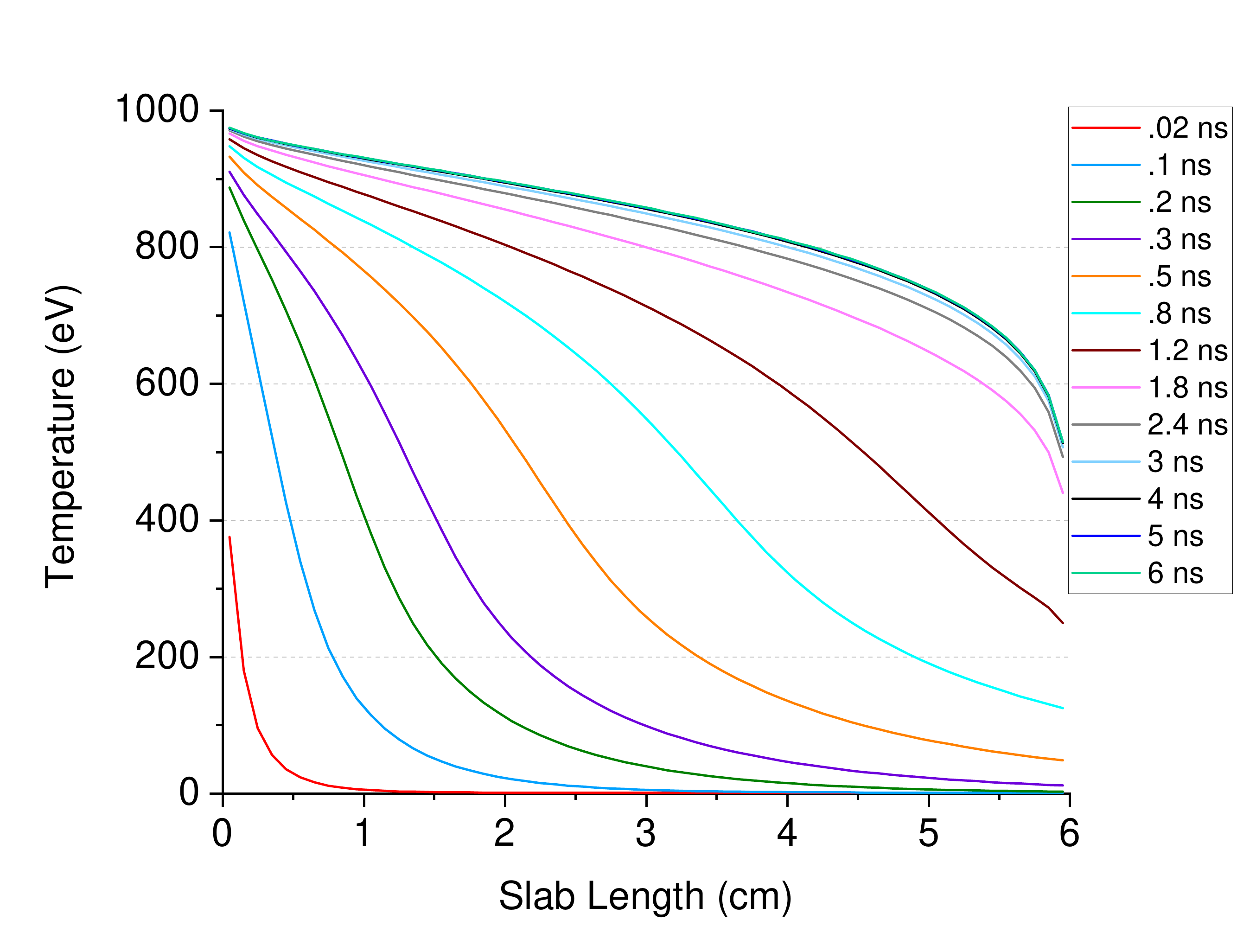}}
	\subfloat[Radiation Energy Density]{\includegraphics[width=.5\textwidth]{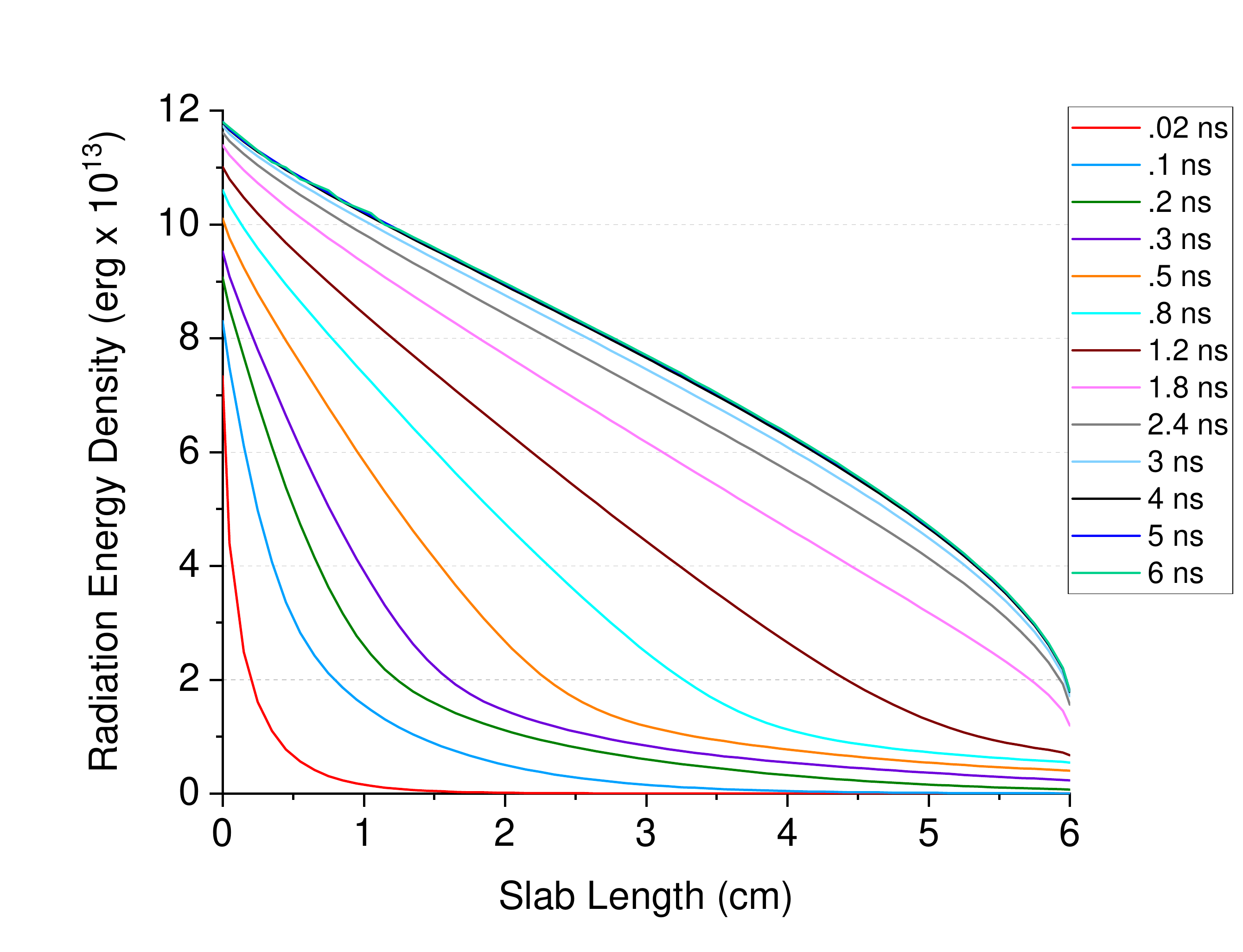}}
	\caption{ \label{fig:sol-eps=1e-5} The solution to the F-C test on the given grid in phase space and time generated by the QD-PODG model with $\varepsilon=10^{-5}$.}
	\vspace*{-0.25cm}
	\subfloat[Temperature]{\includegraphics[width=.5\textwidth]{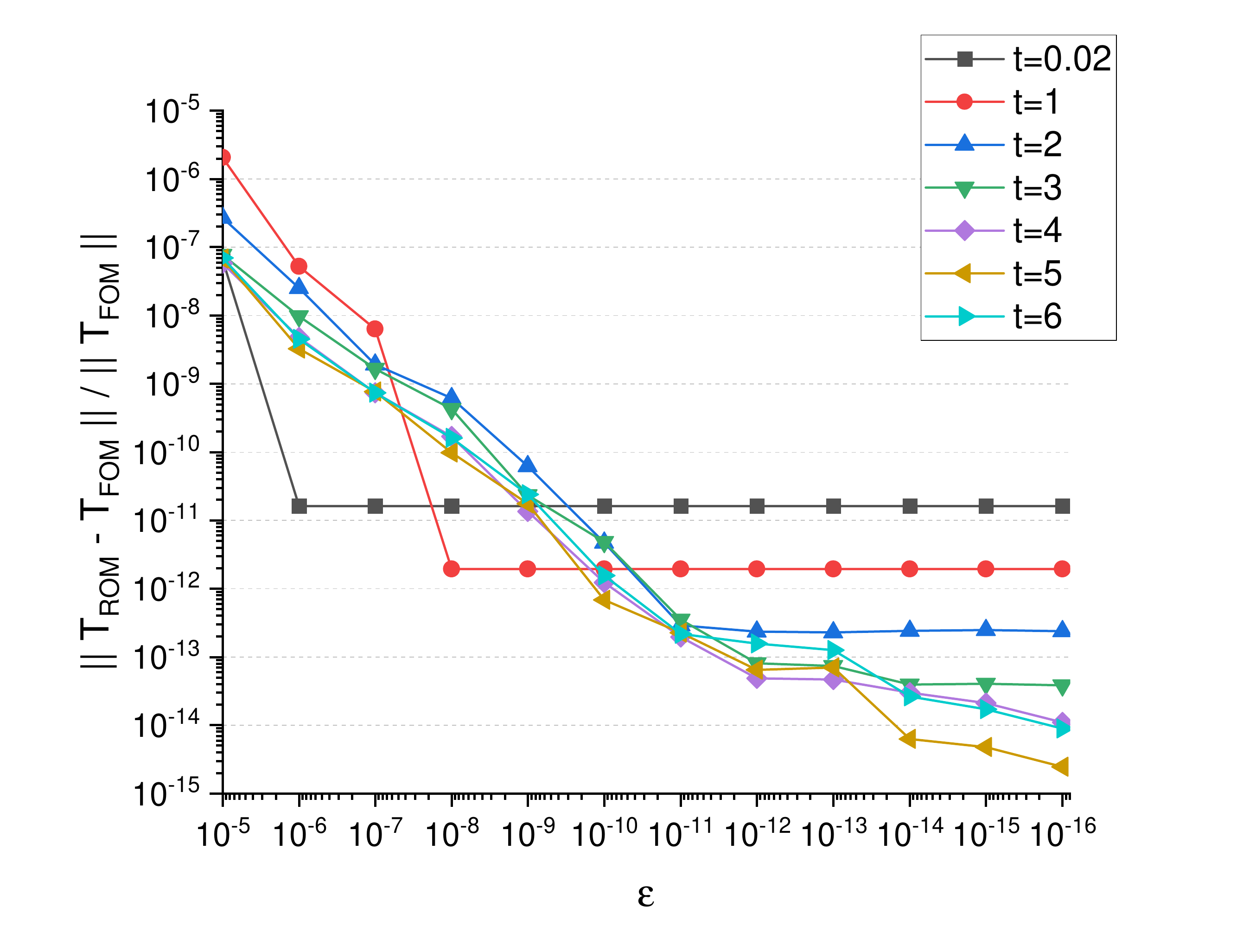}}
	\subfloat[Radiation Energy Density]{\includegraphics[width=.5\textwidth]{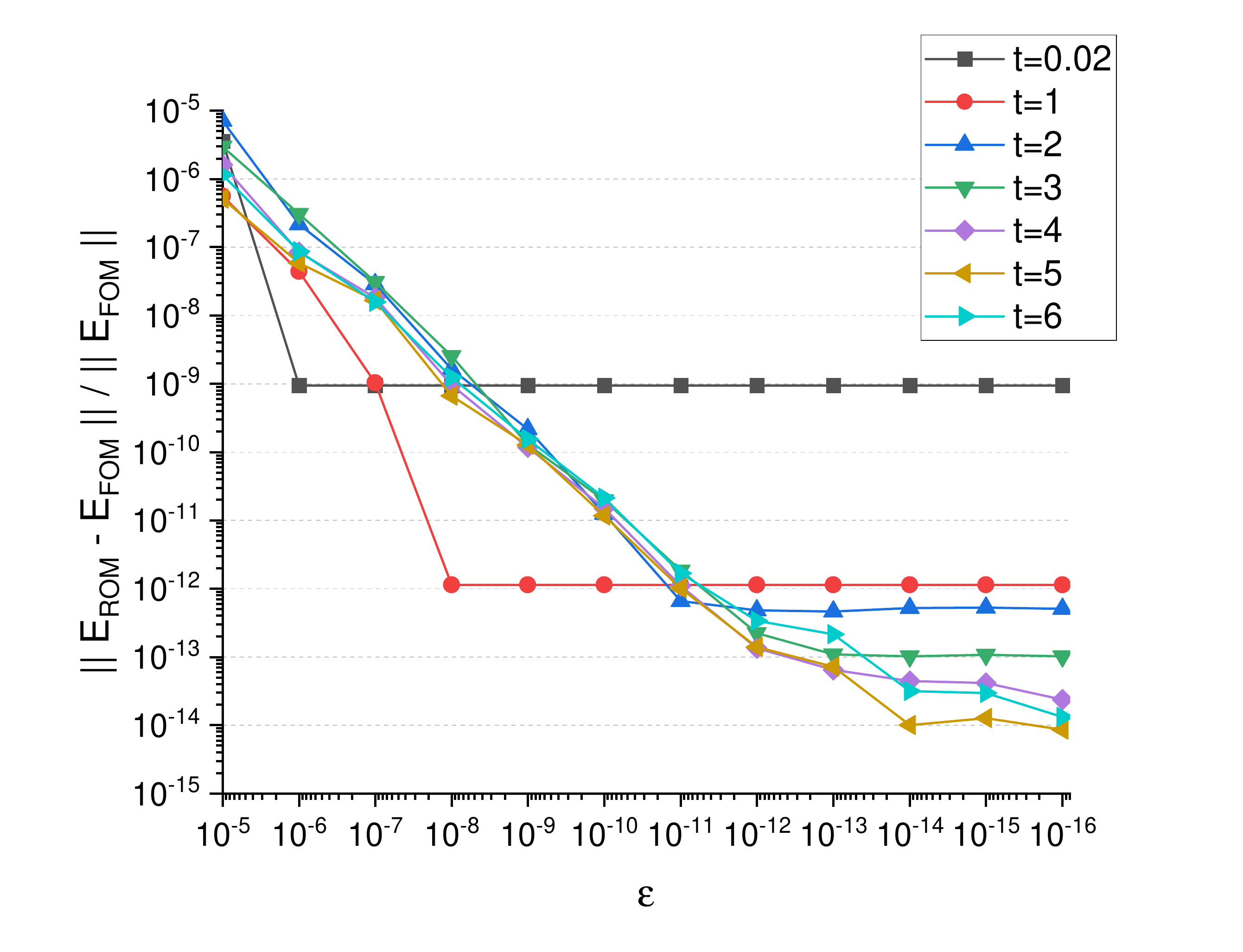}}
	\caption{Relative errors of the QD-PODG model compared to the FOM solution in the 2-norm vs. $\varepsilon$ \label{fig:time_err_trend}}
\end{figure}

Figure \ref{fig:time_err} shows that as $\varepsilon$ decreases, the relative error of the QD-PODG model trends downward as well.
Upon inspection, one can see that the ROM with $\varepsilon=10^{-6},10^{-7}$ is
exceptionally
  accurate for $t\leq 0.3$ ns. This comes from the full-rank basis representation of $\mathbf{A}_1$ that occurs for all $\varepsilon<10^{-5}$, as was shown in Figure \ref{fig:exp_ranks}. Similarly, the high accuracy for $t\leq 1.2$ while using $\varepsilon< 10^{-7}$ follows from the fact that the full-rank basis representation of  $\mathbf{A}_2$ is used for $\varepsilon<10^{-7}$.
Considering overall accuracy, even with very low-rank ($\varepsilon=10^{-5}$) the QD-PODG model maintains a relative error in both material temperature and radiation energy density below $10^{-5}$.
Figure \ref{fig:sol-eps=1e-5} depicts the solution to the F-C test generated with the  QD-PODG model, using the criterion from equation \eqref{eta-g} as $\varepsilon=10^{-5}$.
This ROM can be compared to the relative errors in the 2-norm of the popular multigroup $P_1$ (mg-$P_1$) and multigroup flux-limited diffusion (mg-FLD) ROMs \cite{olson-auer-hall-2000} found for the same test problem, also shown in Figure \ref{fig:time_err}. The results show that even with $\varepsilon=10^{-5}$ the QD-PODG model yields a far more accurate solution than these other ROMs by
 roughly 3-4 orders of magnitude. Let us also take note that when using all POD modes ($\varepsilon=10^{-16}$) the QD-PODG model converges to the  FOM  solution within the iterative convergence bounds with the exception of the radiation energy density while $t<0.5$ ns.
This comes from higher errors found at the radiation wavefront during formation, which is a difficult process to capture given how rapidly it progresses and can be prone to larger numerical errors than other parts of the solution.

Similarly to Figure \ref{fig:time_err}, Figure \ref{fig:time_err_trend} also displays the relative error in the solution of the F-C test obtained by the QD-PODG ROM compared to the FOM solution in the 2-norm, but plotted in a unique format.  In Figure \ref{fig:time_err_trend} each curve corresponds to a specific instant of time, showing how the error of the QD-PODG model changes with respect to $\varepsilon$ when time is held static. This plot clearly demonstrates the convergence behavior of the ROM solution as $\varepsilon$ decreases to zero.

\section{\label{sec:end} Conclusions}

In this paper,  we  presented a new  ROM   for high-energy density TRT problems.
The proposed methodology is based on
 the nonlinear projection approach and Galerkin projection combined with the POD.
The developed ROM efficiently reduces dimensionality of TRT problems and was shown
capable of producing solutions with various levels of fidelity. The accuracy varies based on the rank of the POD basis used to project the BTE, and the  ROM solution converges to the FOM solution as this rank is increased. As such the developed ROMs enable the use of practical and efficient simulations by significantly reducing dimensionality of the problem while maintaining sufficient accuracy. The ROMs presented here also possess the capability for parameterization, which is an avenue the authors will be pursuing in the future.

The promising performance of the QD-PODG model motivates further research on this approach.
An extension to  2D geometry is the next logical step. To make the method robust for such extensions, a desirable feature is to enforce positivity of the expanded intensities. Secondly, work must be done towards generation of
enhanced  POD for the problems at hand; one technique would be to use symmetry-reduction methods \cite{rowley-mardsen-2000,reiss-2004} which are known to improve basis generation for traveling waves.

\section*{Acknowledgments}

This research project  is funded  by the Department of Defense, Defense Threat Reduction Agency, grant number HDTRA1-18-1-0042.
The content of the information does not necessarily reflect the position or the policy of the federal government, and no official endorsement should be inferred.

\bibliography{QD-PODG-TRT-arXiv}

\begin{thebibliography}{10}
\expandafter\ifx\csname url\endcsname\relax
  \def\url#1{\texttt{#1}}\fi
\expandafter\ifx\csname urlprefix\endcsname\relax\def\urlprefix{URL }\fi
\expandafter\ifx\csname href\endcsname\relax
  \def\href#1#2{#2} \def\path#1{#1}\fi

\bibitem{jc-dya-m&c2019}
J.~Coale, D.~Anistratov, A reduced-order model for thermal radiative transfer
  problems based on multilevel quasidiffusion method, in: Int. Conf. on
  Mathematics and Computational Methods Applied to Nuclear Science and
  Engineering (M\&C 2019), Portland, OR, August 25, 2019, pp. 278--287.

\bibitem{jc-dya-tans2019}
J.~Coale, D.~Y. Anistratov, Data-driven grey reduced-order model for thermal
  radiative transfer problems based on low-order quasidiffusion equations and
  proper orthogonal decomposition, Transaction of American Nuclear Society 121
  (2019) 836--839.

\bibitem{ragusa-morel-m&c2019}
P.~A. Behne, J.~C. Ragusa, J.~E. Morel, Model order reduction for {S}$_n$
  radiation transport, in: Int. Conf. on Mathematics and Computational Methods
  Applied to Nuclear Science and Engineering (M\&C 2019), Portland, OR, August
  25, 2019, pp. 2481--2490.

\bibitem{ryan-jcp2020}
Z.~Peng, R.~G. McClarren, M.~Frank, A low-rank method for two-dimensional
  time-dependent radiation transport calculations, Journal of Computational
  Physics 421 (2020) 109735.

\bibitem{Choi-2021}
{Y. Choi}, {P. Brown}, {B. Arrighi}, {R. Anderson}, Space-time reduced order
  model for large-scale linear dynamical systems with application to boltzmann
  transport problems, Journal of Computational Physics 424 (2021) 109845.

\bibitem{sirovich-1987}
{L. Sirovich}, Turbulence and the dynamics of coherent structures. parts i-iii,
  Quarterly of Applied Mathematics XLV (1987) 561--590.

\bibitem{berkooz-holmes-lumley-1993}
{G. Berkooz}, {P. Holmes}, {J. L. Lumley}, The proper orthogonal decomposition
  in the analysis of turbulent flows, Annual Review of Fluid Mechanics 25
  (1993) 539--575.

\bibitem{volkwein-2002}
{K. Kunisch}, {S. Volkwein}, Galerkin proper orthogonal decomposition methods
  for a general equation in fluid dynamics, SIAM J. Numer. Anal 40 (2002)
  492--515.

\bibitem{Benner-siam-2015}
P.~Benner, S.~Gugercin, K.~Wilcox, A survey of projection-based model reduction
  methods for parametric dynamical systems, SIAM Review 57 (2015) 483--531.

\bibitem{mla-ttsp-1997}
M.~L. Adams, Subcell balance methods for radiative transfer on arbitrary grids,
  Transport Theory and Statistical Physics 26 (1997) 385--431.

\bibitem{volkwein-2013}
S.~Volkwein, Model reduction using proper orthogonal decomposition, lecure
  Notes, University of Konstanz (2013).

\bibitem{gol'din-1964}
{V. Ya. Gol'din}, A quasi-diffusion method of solving the kinetic equation,
  USSR Comp.\ Math.\ and Math.\ Phys. 4 (1964) 136--149.

\bibitem{gol'din-1972}
{V. Ya. Gol'din}, {B. N. Chetverushkin}, Methods of solving one-dimensional
  problems of radiation gas dynamics, USSR Comp.\ Math.\ and Math.\ Phys. 12
  (1972) 177--189.

\bibitem{PASE-1986}
{V. Ya. Gol'din}, {D. A. Gol'dina}, {A. V. Kolpakov}, {A. V. Shilkov},
  Mathematical modeling of hydrodynamics processes with high-energy density
  radiation, Problems of Atomic Sci. \& Eng.: Methods and Codes for Numerical
  Solution of Math. Physics Problems 2 (1986) 59--88, in Russian.

\bibitem{auer-mihalas-1970}
L.~H. Auer, D.~Mihalas, On the use of variable {E}ddington factors in non-{LTE}
  stellar atmospheres computations, Monthly Notices of the Royal Astronomical
  Society 149 (1970) 65--74.

\bibitem{dya-jcp-2019}
D.~Y. Anistratov, Stability analysis of a multilevel quasidiffusion method for
  thermal radiative transfer problems, Journal of Computational Physics 376
  (2019) 186--209.

\bibitem{fleck-1971}
{J. A. Fleck}, {J. D. Cummings}, An implicit monte carlo scheme for calculating
  time and frequency dependent nonlinear radiation transport, Journal of
  Computational Physics 8 (1971) 313--342.

\bibitem{rowley-mardsen-2000}
{C. W. Rowley}, {J. E. Mardsen}, Reconstruction equations and the
  karhunen-lo\`eve expansion for systems with symmetry, Physica D 142 (2000)
  1--19.

\bibitem{reiss-2004}
{J. Reiss, P. Schulze, J. Sesterhenn and V. Mehrmann}, The shifted proper
  orthogonal decomposition: A mode decomposition for multiple transport
  phenomena, SIAM Journal of Scientific Computing 40 (2018) A1322--A1344.

\bibitem{olson-auer-hall-2000}
{G. L. Olson, L. H. Auer and M. L. Hall}, Diffusion, {$P_1$}, and other
  approximate forms of radiation transport, Journal of Quantitative
  Spectroscopy \& Radiative Transfer 64 (2000) 619--634.

\end{thebibliography}
\bibliographystyle{elsarticle-num}

\end{document}